\documentclass[11pt]{amsart}
\usepackage{bm}
\usepackage{amssymb}
\usepackage{graphicx}	
\usepackage{xcolor}
\usepackage{amsmath}
\usepackage{enumerate}
\usepackage{blindtext}
\usepackage{multirow}
\usepackage{mathrsfs}

\newtheorem{thm}{Theorem}[section]
\newtheorem{Def}[thm]{Definition}
\newtheorem{rem}[thm]{Remark}
\newtheorem{cor}[thm]{Corollary}

\newtheorem{ex}[thm]{Example}
\newtheorem{prob}[thm]{Problem}
\newtheorem{lemma}[thm]{Lemma}
\newcommand{\bdfn}{\begin{Def} \rm}
\newcommand{\edfn}{\end{Def}}

\newcommand{\tfae}{the following are equivalent}

\newcommand{\ra}{\rightarrow}
\newcommand{\Ra}{\Rightarrow}

\newcommand{\ci}{\subseteq}

\newcommand{\al}{\alpha}

\newcommand{\de}{\delta}

\newcommand{\la}{\lambda}

\newcommand{\La}{\Lambda}

\newcommand{\beqa}{\begin{eqnarray*}}
\newcommand{\eeqa}{\end{eqnarray*}}
\newcommand{\vertiii}[1]{{\left\vert\kern-0.25ex\left\vert\kern-0.25ex\left\vert #1
    \right\vert\kern-0.25ex\right\vert\kern-0.25ex\right\vert}}
\newcounter{cnt1}
\newcounter{cnt2}
\newcounter{cnt3}
\newcounter{cnt4}
\newcommand{\blr}{\begin{list}{$($\roman{cnt1}$)$} {\usecounter{cnt1}
\setlength{\topsep}{0pt} \setlength{\itemsep}{0pt}}}
\newcommand{\blR}{\begin{list}{\Roman{cnt4}.\ } {\usecounter{cnt4}
\setlength{\topsep}{0pt} \setlength{\itemsep}{0pt}}}
\newcommand{\bla}{\begin{list}{$(\alph{cnt2})$} {\usecounter{cnt2}
\setlength{\topsep}{0pt} \setlength{\itemsep}{0pt}}}
\newcommand{\bln}{\begin{list}{$($\arabic{cnt3}$)$} {\usecounter{cnt3}
\setlength{\topsep}{0pt} \setlength{\itemsep}{0pt}}}
\newcommand{\el}{\end{list}}

\begin{document}
\title[ Order-preserving unique extensions]{Order-preserving unique Hahn-Banach extensions}
\author[Paul] {Tanmoy Paul$^{^{*}}$}
{\thanks{$^{*}$Corresponding Author}}
\address[Tanmoy Paul]
{Department of Mathematics\\ I. I. T Hyderabad\\ Telangana\\India,\textit{E-mail~:}\textit{tanmoy@math.iith.ac.in}}
\author[Rao]{T. S. S. R. K. Rao}
\address[T. S. S. R. K. Rao]
{Department of Mathematics\\
Shiv Nadar Institution of Eminence \\
Gautam Buddha Nagar, UP-201314 \\ India,
\textit{E-mail~:}
\textit{srin@fulbrightmail.org}}
\begin{abstract}
	Let $X$ be a real Banach lattice with a unit, and let $Y \subseteq X$ be a closed subspace containing the unit.  In this paper, we study the order-theoretic (also isometric) structure of $Y$ that it may inherit from $X$ under some additional conditions. Specifically, we consider the case where every continuous positive linear functional in the unit sphere of $Y^\ast$ admits a unique positive norm-preserving extension in $X^\ast$. Our results depend on the specific nature of the embedding of $Y$ in $X$. For a compact convex set $K$ with closed extreme boundary $\partial_e K$, for the restriction isometry of $A(K)$ (which is also order-preserving) into $C(\partial_e K)$, uniqueness of extensions of positive functionals leads to $K$ being a simplex and the restriction embedding being onto. In contrast, for a Choquet simplex $K$, we show that under the canonical embedding in the bidual $A(K)^{\ast\ast}$ (which is an abstract $M$-space), positive linear functionals having unique positive linear extensions imply that $K$ is a finite-dimensional simplex. \end{abstract}
\maketitle
\section{Introduction}
Let $K$ be a compact convex set in a locally convex topological vector space and let $A(K)$ denote the space of real-valued affine continuous functions on $K$, equipped with the supremum norm. For a compact space $\Omega$, let $C(\Omega)$ denote the space of real-valued continuous functions on $\Omega$. We will be using basic convexity theory from
\cite{A} Chapter I and  \cite{AE} Chapter II, Section 6. Some of the essential facts we will use are recorded in Section 2.  We refer to the monograph \cite{L}, Chapter I, for results from order structure theory, geometric properties, and classification of abstract $L$ and $M$-spaces. See also the monograph \cite{Lu}.
\vskip 1em
 We recall from Chapter 1 of \cite{L} that the bidual of an abstract $M$-space $X$  is always a $C(\Omega)$-space with the constant function ${\bf 1}$ as the unit.
Thus, a space $X$ and its bidual have lattice structures, and the canonical embedding of $X$ into its bidual is an order-preserving map. An interesting open problem is whether $X$ inherits the geometric structure of a $C(\Omega)$-space under the assumption of unique order-preserving and norm-preserving extensions of linear functionals from $X$ to its bidual?
\vskip 1em
In this article, we attempt to solve this for order unit spaces, with an order unit. See the recent paper \cite{KS} for possible approaches when there is no order unit. The article \cite{KRS} has developed new classes of function spaces that are similar to $A(K)$-spaces. It is possible that the integral representation approach that we use here, via Choquet's theorem, can also be applied to these classes.  See also the Zb Math Open reviews of these articles by Rao.

\vskip 1em
Let $X$ be a real Banach space and let $Y \subseteq X$ be a closed subspace. Suppose $Y$ is a $U$-subspace of $X$, in the sense that non-zero elements of $Y^\ast$ have unique norm-preserving extensions in $X^\ast$ (see \cite{P}). An interesting question is to consider what additional geometric conditions need to be assumed on $Y$ so that some of the geometric properties of $X$ are inherited by a $U$-subspace $Y$? A weaker condition than that of a $U$-subspace has also been studied in the literature. Here, one seeks the uniqueness of a norm-preserving extension for only norm-attaining functionals in $Y^\ast$, see \cite{Lim}. In such a case, we say that $Y$ is a $wU$-subspace of $X$. See the recent paper  \cite{DPR} for stability results for $U$ and $wU$ subspaces. We note that the unique extensions we study in this paper are also norm-attaining functionals. Section 3 considers some implications of these geometric properties, with particular emphasis on spaces $X$ that contain an $A(K)$ space as a $wU$ subspace (Corollary 3.5).
\vskip 1em
See also \cite{KPR}, where we have considered the uniqueness of extensions for algebraic embedding of  $c_0$ in the space of compact operators ${\mathcal K}(\ell^p)$, $1 \leq p < \infty$, $p \neq 2$, analogous to the multiplication operator embedding of $c_0$, in the case of Hilbert spaces.
\vskip 1em
We recall the following from Chapter II of \cite{A}.
\begin{Def}
\bla
\item A partially ordered vector space $A$ over $\mathbb{R}$ is said to be {\bf Archimedean} if the negative elements $a\in A^-$ are the only ones for which $\{\alpha a:\alpha\in\mathbb{R}^+\}$ has an upper bound.
\item A positive element $e\in A$ is said to be an {\bf order unit} if the smallest order ideal generated by $e$ is $A$.
\el
\end{Def}
Let $X$ be a norm-complete order unit Archimedean space. By Kadison's theorem (\cite{AE} Theorem II.1.8), $X$ is order and linearly isometric to an $A(K)$-space for a compact convex set $K$. We recall that any $C(\Omega)$-space is an $A(K)$-space, where $K$ is identified as the set of probability measures in $C(\Omega)^\ast$, with the weak$^\ast$-topology. The results in Section 4 illustrate how the structure of a compact convex set can be determined by considering unique extensions of positive functionals in $A(K)^\ast$  to positive functionals on a $C(\Omega)$-space. The embedding of $A(K)$ into $C(\Omega)$ or to the bidual $A(K)^{\ast\ast}$ under the canonical embedding turns out to be the determining factor.
\vskip 1em
   Let $K$ be a compact convex set. Suppose $K$ is not a singleton, then for any $k_1 \neq k_2$ in $K$, the discrete measures $\delta_{\frac{k_1+k_2}{2}}$, $\frac{\delta_{k_1}+\delta_{k_2}}{2}$ are two distinct positive extensions in $C(K)^*$ of the evaluation map $\delta_{\frac{k_1+k_2}{2}}$ on $A(K)$. Thus, it is important to choose carefully an embedding of $A(K)$ in a $C(\Omega)$ space to determine the structure of $K$ using the uniqueness of positive extensions.
  \vskip 1em
    Let $\Phi: A(K) \rightarrow C(\Omega)$ be an order preserving linear isometry such that $\Phi({\bf 1})= {\bf 1}$. In Theorem~\ref{T9}, we show that if positive (and hence continuous) linear functionals on $\Phi(A(K))$ have unique extensions to positive functionals in $C(\Omega)^\ast$, then $K$ is a Bauer simplex, i.e., a Choquet simplex whose extreme boundary $\partial_e K$ is a closed set.
  \vskip 1em
 A particular embedding that has been well studied from the point of view of unique extensions is the canonical embedding of a Banach space $X$ in its bidual $X^{\ast\ast} $. In this embedding, if $X$ is a $U$-subspace of $X^{\ast\ast}$, then every non-zero functional $x^\ast \in X^\ast$ has unique norm preserving extension $x^\ast \in X^{\ast\ast\ast}$ (here again we are considering the canonical embedding of $X^\ast \subseteq X^{\ast\ast\ast})$. When $X$ has this property, it is called a Hahn-Banach smooth space, see \cite{S}.
    \vskip 1em
   For a Banach space $X$, let $X_1$ and $S(X)$ denote the closed unit ball and the unit sphere, respectively. It follows from \cite{G}  (see also Lemma III.2.14 in \cite{HWW}) that  Hahn-Banach smoothness of $X$ is equivalent to the identity map $i: (S(X^\ast),weak^\ast) \rightarrow (S(X^\ast),weak)$ being continuous. We further note that as an application of $weak^*$ lower-semi-continuity of the norm, for a non-reflexive Banach space $X$, this is also the same as continuity of the map, $i: (X^\ast_1,weak^\ast) \rightarrow (X^\ast_1,weak)$, see \cite{DMR}. Chapter III of \cite{HWW}   contains several examples of spaces of operators and function spaces that are Hahn-Banach smooth spaces. In Section 3, we explore geometric implications of the identity map on $S(X^\ast)$ having points of weak$^\ast$-weak continuity.  Since a normed linear space and its completion have isometric dual spaces, and since the weak$^\ast$-topologies on the dual unit balls of the normed linear space and its completion agree, some of the results in Section 3 are also valid for normed linear spaces.
     \vskip 1em
  In Section 4, for a compact Choquet simplex $K$, we consider the canonical embedding of $A(K)$ in its bidual. We recall from Section 9,  Theorem 2 in \cite{AE} that in this case $A(K)^{\ast\ast}$ is a lattice. The question, when do positive norm-attaining linear functionals have unique extensions to positive functionals on the bidual, gets answered in Theorem 4.1, when we show that this happens only when the canonical embedding is onto. Since $K$ is a simplex, $A(K)^\ast$ is isometric to $L^1(\mu)$ for some positive measure $\mu$. See Section 18 in \cite{L}. Thus, when $A(K)$ is reflexive, it is isometric to the finite-dimensional space, $\ell^\infty(k)$ for some positive integer $k$.

  \section{Preliminaries}
  In this section, we recall some preliminaries from the theory of compact convex sets that we will be using and motivate the results in the subsequent sections. These can be found in \cite{AE} and \cite{A}. Let $K$ be a compact convex set in a locally convex space $E$.
  Let $S=\{\Lambda \in A(K)^\ast_1:\Lambda(1)=1\}$ be equipped with the weak$^\ast$-topology. This is called the state space of $A(K)$ and it is a face, i.e.,  an extreme, convex subset of $A(K)^\ast_1$. Consider $A(K) \subseteq C(K)$. For $k \in K$, let $\delta_k$ denote the evaluation functional, which will be interchangeably considered as a Dirac measure on $K$.
  We recall that any norm preserving extension of  $ \delta_k \in S$ to $C(K)$ is a regular Borel probability measure $\mu$ on $K$ and $a(k) = \int_K a~ d\mu$ for all $a \in A(K)$. We write this as the resultant equation, $r(\mu) = k$.  Then it follows from Proposition I.2.1 in \cite{A} that $k \mapsto \delta_k$ is a homeomorphism of $K$ onto $S$. It is also an affine map, since for any $k_1,k_2 \in K$, $\lambda \in [0,1]$, $(\lambda \delta_{k_1} + (1-\lambda) \delta_{k_2})(a) = a(\lambda k_1 + (1-\lambda) k_2)= \lambda a(k_1)+ (1-\lambda) a(k_2)$, for all $a \in A(K)$. So we will just write $K$ instead of $S$ while working in $A(K)^\ast$.
\vskip 1em
For a compact convex set $K$, let $\partial_e K$ denote the set of extreme points of $K$. We recall that as a consequence of the Krein-Milman theorem, we have $\|a\|=\sup\{|a(k)|: k \in K\}=\sup\{|a(k)|:k \in \partial_e K\}$. If $\partial_e K$ is a closed set, the mapping $a \mapsto a|\partial_e K$ is called the restriction embedding of $A(K)$ in $C(\partial_e K)$. We always consider the dual unit ball, $A(K)^\ast_1$, equipped with the weak$^\ast$-topology, and closures of subsets here are always with respect to the weak$^\ast$-topology.
\vskip 1em
If $D$ is a convex set and $C \subset D$, we denote by $co(C)$ the smallest convex subset of $D$ containing $C$. Using the Jordan decomposition of measures in $C(K)^\ast$ as the difference of positive and negative parts, as in the proof of Proposition I.2.1 of \cite{A}, we see that $A(K)^\ast_1 = co(K \cup -K)$ (convex hull).
\vskip 1em
 Thus  $\partial_e A(K)^\ast_1 = \partial_e K \cup - \partial_e K$. In particular we note that $\overline{\partial_e A(K)^\ast_1} \subseteq K \cup -K$ (closure taken in the weak$^\ast$-topology) and hence $\overline{\partial_e A(K)^\ast_1 } \subseteq S(A(K)^\ast)$ and elements of the set $\overline{\partial_e A(K)^\ast_1 }$ are all norm attaining functionals.
\vskip 1em
Let ${\bf 1}$ denote the constant function. Since $a \in A(K)$ is in $A(K)_1$ if and only if in the point-wise ordering of $A(K)$,  ${\bf -1 } \leq  a  \leq {\bf 1}$ we see that any $A(K)$ space under the point-wise ordering is an order unit space with the constant function ${\bf 1}$ as the order unit (see \cite {A}, Chapter II).
When $K$ is a Choquet simplex, $A(K)$ has the Riesz-decomposition property and $A(K)^\ast$ is an abstract $L$-space (see Theorem III.7.1 in \cite{AE} or Section 19, Theorem 2 in \cite{L}). The monograph \cite{AE} provides several examples of complex function spaces whose state spaces are Choquet simplexes.
  \vskip 1em
  We also need the notion of a parallel face from \cite{AE}.
\bdfn
Let $K$ be a compact convex set in a locally convex space $E$. Let $F$ be a closed face of $K$, i.e., $F$ is a convex set and if $\lambda k_1+(1-\lambda)k_2 \in F$ for $k_1,k_2 \in K$ and $\lambda \in (0,1)$, then $k_1, k_2 \in F$. The complementary set $F'$ is the union of faces disjoint from $F$. Suppose $F'$ is also a convex set. $F$ is said to be a {\bf  parallel face}, if for each $k\in K$ there exists $x\in F$, $y\in F^\prime$ and a unique $0\leq\la\leq 1$ such that $k=\la x+(1-\la)y$. If $x,y$ are also unique in the above decomposition, then $F$ is said to be a split face.
\edfn
It is known that for finitely many closed split faces
$\{F_i\}^n_{i=1}$ of $K$, the convex hull $co(\bigcup_{i=1}^n F_i)$ is a closed split face, see Corollary II.6.8 in \cite{A}.
\vskip 1em
Let $\Omega$ be a compact set. A motivating example for this investigation is the set $P(\Omega)$ of probability measures on $\Omega$, with the weak$^\ast$-topology. This is a Choquet simplex with extreme boundary identified with $\Omega$. Any closed face $F$ is of the form $\overline{co}(E)$ for a closed set $E \subseteq \Omega$. In this case, $F'$ is the set of probability measures concentrated on $\Omega \setminus E$. $F'$ is a face and $F$ is a split face of $P(\Omega)$.
\section{Geometric implications of unique extensions}
We first study the nature of unique Hahn-Banach extensions depending on the embedding of the subspace $Y$ into $X$. These results motivate some of the conditions we will assume for ordered spaces. See also the prelude in Section 4.
\vskip 1em
We recall from the introduction that  $x^\ast \in S(X^\ast)$ has a unique norm preserving extension in $S(X^{\ast\ast\ast})$ if and only if  $x^\ast$ is a point of weak$^\ast$-weak continuity for the identity mapping on $S(X^\ast)$. Our first result shows that both weakly Hahn-Banach smoothness and Hahn-Banach smoothness are hereditary properties. For a closed subspace $Y \subseteq X$, in what follows, we use the canonical embeddings and identifications of $Y^{\ast\ast} = Y^{\bot\bot} \subseteq X^{\ast\ast}$. Since the weak$^\ast$-topology is, in general, not metrizable, we will use the notions of nets and subnets from \cite{K} and, for notational simplicity, write only the indexing set of the net, say $\Delta$. Also, since we most often work in compact, Hausdorff spaces, rather than going through subnets, we assume, without loss of generality, that the nets under consideration converge.  We also recall that for Hausdorff spaces, $\Omega_1,~\Omega_2$, a function $f: \Omega_1 \rightarrow \Omega_2$ is continuous if and only if any net $\{\omega_{\alpha}\}_{\alpha \in \Delta}  \subset \Omega_1$ such that $\omega_{\alpha} \rightarrow \omega$ has a subnet ( still indexed by $\Delta$),  $\{\omega_{\alpha'}\}_{\alpha' \in \Delta}$ such that $f(\omega'_{\alpha}) \rightarrow f(\omega)$.
\begin{thm}\label{T1} Let $Y \subseteq X$ be a closed subspace of a Banach space $X$. Suppose $i: (S(X^\ast),weak^\ast) \rightarrow (S(X^\ast), weak)$ is continuous. Then the same conclusion holds for the identity map on $S(Y^\ast)$.	
\end{thm}
\begin{proof}
Let $y^\ast_0 \in S(Y^\ast)$. To show that it is a point of continuity, let $\{y^\ast_{\alpha}\}_{\alpha \in \Delta} \subseteq S(Y^\ast)$ be a net such that $y^\ast_{\alpha} \rightarrow y^\ast_0$ in the weak$^\ast$-topology of $Y^\ast$.  To prove weak continuity, it is enough to exhibit a subnet of the given net converging weakly to $y^\ast_0$.
\vskip 1em
For each $\alpha$, choose a norm-preserving Hahn-Banach extension  $x^\ast_{\alpha}\in S(X^*)$ of $y^\ast_{\alpha}$. By going through a subnet if necessary, we may and do assume that $x^\ast_{\alpha} \rightarrow x^\ast$ in the weak$^\ast$-topology of $X^\ast_1$.  As $x^\ast= y^\ast_0$ on $Y$, we have, $\|x^\ast\|= 1 $. Thus by hypothesis, $x^\ast_{\alpha} \rightarrow x^\ast$ in the weak topology of $X^\ast$.
\vskip 1em
Let $\tau \in Y^{\ast\ast}=Y^{\bot\bot} = (X^{\ast}/Y^\bot)^\ast$ and let $\pi: X^\ast \rightarrow X^\ast / Y^{\bot}$ be the quotient map. We have the identification $X^\ast/Y^{\bot} = Y^\ast$.  We note that $\pi(x^\ast_{\alpha}) \rightarrow \pi(x^\ast)$ in the weak topology. Also $\pi(x^\ast_{\alpha})=y^\ast_{\alpha}$ for all $\alpha \in \Delta$, $\pi(x^\ast)=y_0^\ast$. Thus as $\tau \in Y^{\bot\bot} \subseteq X^{\ast\ast}$, using the hypothesis, we see that $\tau(y^\ast_{\alpha}) \rightarrow \tau(y_0^\ast)$. Thus, we have $y^\ast_{\alpha} \rightarrow y^\ast_0$ in the weak topology.
\end{proof}
In the next set of results, we will be using the following application of the Krein-Milman theorem several times. We recall the notion of a face from Section 2.
Let $Y\subseteq X$ be a subspace. The following is a well-known fact; we include it here for completeness.
\begin{lemma}\label{L1}
Let $Y$ be a subspace of a Banach space $X$. Then any $y^\ast \in \partial_e Y^\ast_1$ has an extension to $x^\ast \in \partial_e X^\ast_1$.
\end{lemma}
\begin{proof}
The set of norm-preserving extensions of $y^\ast\in\partial_e Y^\ast_1$ is a weak$^\ast$-closed convex subset of $X_1^*$. Also if $x_i^*\in X_1^*$, $(i=1,2)$ are any two functionals satisfying ${\lambda x_1^*+(1-\lambda)x_2^*}|_Y=y^*$, as $y^\ast$ is an extreme point, $x_i^*$'s are also norm preserving Hahn-Banach extensions of $y^*$. Thus, the set of norm-preserving extensions is a weak$^\ast$-closed face of $X_1^*$.
Hence by the Krein-Milman theorem, the set of Hahn-Banach extensions of $y^*$ contains an extreme point of $X^\ast_1$.
\end{proof}
\begin{ex}
	Even when $Y \subseteq X$ is a $U$-subspace, the unique extension of a point of weak$^\ast$-weak continuity in $S(Y^\ast)$ need not be a point of weak$^\ast$-weak continuity for the identity map on $S(X^\ast)$. We recall that $x_0 \in S(X)$ is a smooth point if there is a unique functional $x^\ast_0 \in \partial_e X^\ast_1$ such that $x_0^\ast(x_0)=1$.  Let $X= C([0,1])$, for $s \in [0,1]$, let $\delta_s$ denote the Dirac measure at $s$. It is well known that $\partial_e C([0,1])^\ast_1 = \{\alpha \delta_s: s\in [0,1],~|\alpha|=1\}$.  We note that $g \in S(C([0,1]))$ is a smooth point if and only if there is a unique $t \in [0,1]$ such that $|g(t)|=1 $.  Since $[0,1]$ has no isolated points, we note that no element of $\partial_e C([0,1])^\ast_1$ is a point of weak$^\ast$-weak continuity for the identity map on $S(C([0,1])^\ast)$. To see this, let $s \in [0,1]$ and choose a sequence of distinct terms different from $s$,  $\{s_n\}_{n \geq 1} \subseteq [0,1]$ such that $s_n \rightarrow s$. We have $\delta_{s_n} \rightarrow \delta_s$ in the weak$^\ast$-topology. If $\delta_s$ is a point of continuity, then $\delta_{s_n} \rightarrow \delta_s$ weakly in $C([0,1])^\ast$. Since $C([0,1])^{\ast\ast}$ contains all bounded, Borel measurable functions on $[0,1]$, by evaluating the Dirac measures at the indicator function, $\chi_{\{s\}}$ , as $s_n \neq s$, we get a contradiction.
\vskip 1em
 For any $f \in C([0,1])$ with $0 \leq f \leq 1$ and $f^{-1}(1)=\{1\}$, since $f$ is a smooth point of $C([0,1])$, clearly the one dimensional space, $Y = span\{f\}$ is a $U$-subspace of $C([0,1])$.
\end{ex}
Let $Y \subseteq X$ be a $U$-subspace. Our next theorem illustrates what we mean by `$Y$ may inherit the geometric structure of $X$'. More precisely Theorem~\ref{T3} ensures if $\partial_e X_1^\ast$ is weak$^\ast$-closed in $S(X^*)$ then the family of functionals in $\partial_e Y_1^\ast$ which have unique Hahn-Banach extensions to $X$ is also relatively weak$^\ast$-closed in $S(Y^\ast)$. Both  Corollaries~\ref{C1} and  \ref{C2} make this situation more evident.
We note that in a dual space, the closure operation is, in general, with respect to the weak$^\ast$-topology. In the proof of the result below, we use the standard fact that in a compact space, if all convergent subnets of a given net converge to the same vector, then the net itself converges to that vector.
\begin{thm}\label{T3}
Let $X$ be a Banach space such that $\partial_e X^\ast_1$ is a weak$^\ast$-closed set. Let $Y \subseteq X$ be a closed subspace and $y^\ast \in (S(Y^\ast) \cap \overline{\partial_e Y^\ast_1})$. Suppose $y^\ast$ has a unique norm-preserving extension to $ x^\ast \in S(X^\ast)$. Then $y^\ast \in \partial_e Y^\ast_1$.
	
\end{thm}
\begin{proof}
Let $\{y^\ast_{\alpha}\}_{\alpha \in \Delta} \subseteq \partial_e Y^\ast_1$ and $y^\ast_{\alpha} \rightarrow y^\ast$ in the weak$^\ast$-topology of $Y^\ast$. Let $\{x^\ast_{\alpha}\}_{\alpha \in \Delta} \subseteq \partial_e X^\ast_1$ be a net of extreme extensions as stated in Lemma~\ref{L1}. As the limit of any weak$^\ast$-convergent subnet is a norm-preserving extension of $y^\ast$, using the uniqueness of the extension $x^\ast$, we see that all weak$^\ast$ convergent subnets of the above net converge to $x^\ast$. Therefore $x^\ast_{\alpha} \rightarrow x^\ast$ in the weak$^\ast$ topology of $X^\ast$. Now by hypothesis, $x^\ast \in \partial_e X^\ast_1$.
	\vskip 1em
Suppose $y^\ast = \frac{y^\ast_1+ y^\ast_2}{2}$ for some $y^\ast_i \in Y^\ast_1$. Since $\|y^*\|=1$, we get $\|y_1^*\|=1=\|y^\ast_2\|$. If $x^\ast_i$ denotes a norm preserving extension in $X_1^\ast$ of the $y_i^\ast$'s, then, as $\frac{x^\ast_1 + x^\ast_2}{2}$ is in the unit ball and is an extension of the unit vector $y^\ast$, it is a norm preserving extension of $y^\ast$, we get $x^\ast = \frac{x^\ast_1 + x^\ast_2}{2}$. Since $x^\ast$ is an extreme point, $x^\ast=x^\ast_1=x^\ast_2$. Thus $y^\ast \in \partial_e Y^\ast_1$.
\end{proof}
Given that $K$ is a face of $A(K)^\ast_1$ (under the canonical embedding), the following corollary can be easily deduced. We again recall that these spaces are considered over the real scalar field. Therefore, as previously mentioned, $\overline {\partial_e A(K)^\ast_1 } \subseteq K \cup - K$, as $K \cup -K$ is weak$^\ast$-closed.
\begin{cor}\label{C1}
Let $K$ be a compact convex set and let $X$ be a Banach space such that $\partial_e X^\ast_1$ is a weak$^\ast$ closed set. If $A(K)$ is isometric to a $wU$-subspace of $X$, then $\partial_e K$ is a closed set.
\end{cor}
\begin{proof}
Ignoring the isometry, without loss of generality, we assume that $A(K) \subseteq X$ is a $wU$-subspace. Let $\omega\in \overline{\partial_e K}$. We identify $\omega$ with the evaluation functional on $A(K)$. From our assumption, it follows that $\omega$ has a unique extension to $X$. From Theorem~\ref{T3} we get, $\omega\in \partial_e K$ and hence the result follows.
\end{proof}
A careful examination of the proof of Theorem 3.4 shows that when $y^\ast$ is also a norm attaining functional, the extensions and the components in the averaging arguments are all norm attaining functionals.
\begin{cor}\label{C2}
Let $X, Y$ be as in Theorem 3.4. Suppose $Y$ is a $wU$-subspace of $X$. Then $\partial_e Y^\ast_1$ contains its norm attaining weak$^\ast$-accumulation points. If $Y$ is such that $\overline{\partial_e Y_1^\ast} \subseteq S(Y^\ast)$ and $Y$ is a $U$-subspace of $X$, then $\partial_e Y^\ast_1$ is a weak$^\ast$-compact set.
\end{cor}
\begin{proof}
Let $y^*$ be a weak$^\ast$-accumulation point of $\partial_e Y_1^*$ which is also norm attaining. Let $\{(y_\alpha^\ast)\}_{\alpha \in \Delta} \subseteq \partial_e Y_1^\ast$ be a net such that $y_\alpha^\ast\rightarrow y^\ast$ in the weak$^\ast$-topology. Let $x_\alpha^\ast\in\partial_e X_1^*$ and $x^\ast\in X^*$ be the norm-preserving extensions of $y_\alpha^*$ and $y^*$, respectively. Since $\partial_e X_1^*$ is weak$^*$-closed, we get $x^*\in \partial_e X_1^\ast$ and this implies $y^\ast\in S(Y^*)\cap \overline{\partial_e Y_1^\ast}$. The first part now follows from Theorem~\ref{T3}. Now suppose that $Y$ is a $U$-subspace of $X$. Since every $y^\ast\in Y^*$ has a unique norm-preserving extension over $X$, the second part is obvious.
\end{proof}

\section{Positive unique extensions}	
 In the context of $A(K)$-spaces we first consider Choquet simplexes $K$, where it is known that $A(K)^{\ast\ast}$ is an abstract $M$-space with an order that is, in the canonical embedding, an extension of the order on $A(K)$ and ${\bf 1}$ continues to be the order unit in $A(K)^{\ast\ast}$. See Chapter 1 in \cite{L}. This raises the question of whether it is possible to classify simplexes by the uniqueness of positive linear extensions from $A(K)$ to its bidual.
\vskip 1em

\begin{thm}\label{T7}
Let $K$ be a Choquet simplex. If positive linear functionals in $A(K)^\ast$ have unique positive linear extensions to the bidual, then $K$ is a finite-dimensional simplex. In particular, when $A(K)$ is an infinite-dimensional space, it is not weakly Hahn-Banach smooth.
\end{thm}
\begin{proof}
Suppose $A(K)$ is an infinite-dimensional space. Since $K$ is a Choquet simplex,  $A(K)^{**}$ is an abstract $M$-space with the constant function $1$ as the unit. As recalled above, the canonical embedding of $A(K)$ in its bidual is an order preserving map (see \cite{L}, Chapter 1 and Sections 19 and 20). Hence, using the identification of $A(K)^\ast$ as boundary measures on $K$ (see the discussion on pg.106 of \cite{A}), one has a projection $P: C(K)^\ast \rightarrow C(K)^\ast$ of norm one such that $\ker(P) = A(K)^\bot$. This projection is the identity mapping on the Dirac measures associated with $\partial_e K$, or equivalently, the evaluation map on $A(K)$. The norm of $A(K)$ is determined by these functionals. It now follows from \cite[Lemma~1(i), Remark~4]{Rao} that $A(K)\ci C(K)\ci A(K)^{**}$, under the canonical embedding of $A(K)$ in its bidual. Therefore, by the hypothesis, $A(K) \subseteq C(K)$, possesses the unique extension property for non-negative functionals by the explanation given at the beginning of this proof. As previously mentioned, for distinct points $k_1, k_2 \in K$, the measure $\de_{\frac{k_1+k_2}{2}} \in A(K)^*$ possesses distinct positive extensions, namely $\de_{\frac{k_1+k_2}{2}}$ and $\frac{\de_{k_1}+\de_{k_2}}{2}$, to the space $C(K)$. This contradiction shows that $A(K)^{\ast\ast} = A(K)$ and hence $A(K)$ is a reflexive space. As noted before, since $A(K)^\ast = L^1(\mu)$ and is also a reflexive space, we conclude that $A(K)^\ast = \ell^1(n)$ for a positive integer $n$. Therefore, $K$ is a finite-dimensional simplex.
\end{proof}

\begin{rem}\label{R1}
Let $K$ be a compact convex set. If $\La\in \partial_e  A(K)^*_1$ is a positive functional, then by Bauer's theorem (Theorem I.6.3 in \cite{AE}), there exists unique representing measure $\mu\in C(K)^\ast$ and $k_0\in \partial_e K$ such that for all $f\in A(K)$, $\La (f)=\mu (f)=f(k_0)$. Clearly in this case $\mu=\de_{k_0}$ is the unique positive Hahn-Banach extension of $\La$ to $C(K)$. Thus, assuming unique extensions for positive extreme functionals is not enough to derive Theorem 4.1.
\end{rem}

The following theorem may be viewed as a new variation on this theme by assuming that the identity map on $S(X^\ast)$ has a weak dense set of points of weak$^\ast$-weak continuity. This is a much weaker assumption since we do not, in general, know the structure of a point of weak$^\ast$-weak continuity.  We note that for closed subspaces $M, N \subset X$, $X= M \bigoplus_1 N$ denotes the $\ell^1$-direct sum, i.e., $\|m+n\|=\|m\|+\|n\|$. $\delta: K \rightarrow A(K)^\ast_1$ denotes the canonical evaluation map and $\delta_k$ its value at $k \in K$.
\begin{thm}\label{T2}
Let $K$ be a compact convex set such that every $k \in \partial_e K$ is a split face of $K$. Suppose the identity map on $S(A(K)^\ast)$ has a weak dense set of points of weak$^\ast$-weak continuity. Then $K$ is a simplex. \end{thm}
\begin{proof}
	It is easy to see that for any $k \in \partial_e K$ by the split face assumption, $A(K)^\ast = span\{\delta_k\} \bigoplus_1 N$  where $N$ is a closed subspace. See page 5, Example I.1.4 (c) in \cite{HWW}. Now for any extreme point $k' \neq k$, as $\delta_{k'} \in N$ we see that for any scalars $\alpha,~\beta$, $
    \|\alpha \delta_k+\beta \delta_{k'}\|= |\alpha|+|\beta|$.  A similar decomposition also works for any finite set $\{\delta_{k_i}\}_{i=1}^n $ of distinct extreme points and their span, so that $\|\sum_{i=1} ^n\alpha_i \delta_{k_i}\|= \sum_{i=1}^n|\alpha_i|$, for any scalars $\alpha_1,...,\alpha_n$. Let $\tau \in S(A(K)^\ast)$ be a point of weak$^\ast$-weak continuity, which we may assume is a point of continuity for the identity map on the dual unit ball (see the remarks from Section 1). Since $A(K)^\ast_1 = \overline {co(\pm\delta((\partial_e K)))}^{weak^\ast}$ and as weak and norm closures of convex sets are the same, $\tau \in \overline {co(\pm (\delta(\partial_e K)))}^{\|.\|}$. Therefore, as points of continuity are weakly dense in $S(A(K)^\ast)$, again using the fact that the weak and norm closures are the same for a convex hull,  we get $A(K)^\ast_1 = \overline{ co(\pm\delta(\partial_e K)))}^{\|.\|}$.
	\vskip 1em
	Consider  $\Gamma = \delta(\partial_e K)$ as a discrete set. We recall that $\ell^1(\Gamma) = \{\alpha \in {\bf R}^{\Gamma} : \sum_{\gamma \in \Gamma}|\alpha(\gamma)|<\infty\}$. It is easy to see that by the assumption on the extreme points of $K$, the canonical map $\Phi: \ell^1(\Gamma) \rightarrow A(K)^\ast$ defined by $\Phi(\alpha)= \sum_{k \in \partial_e K} \alpha(\delta_k)\delta_{k}$ (here we are interpreting coordinates of $\alpha$ as extreme points) is a surjective isometry. Thus, by Theorem 2 in Section 19 from \cite{L}, which states that if  $A(K)^\ast$  is isometric to an $L^1(\mu)$-space, then $K$ is a simplex, we get the conclusion.
\end{proof}

\begin{cor}
	Let $K$ be a compact convex set such that every $k \in \partial_e K$ is a split face of $K$. Suppose $A(K)$ is a weakly Hahn-Banach smooth space. Then $K$ is a finite-dimensional simplex.
\end{cor}
\begin{proof}
	The hypothesis of weakly Hahn-Banach smoothness implies that every point of $K $ is a point of weak$^\ast$-weak continuity. Here we assume that $K$ is identified with the evaluational functional $\{\de_{k}:k\in K\}$. Thus, as before, $A(K)^\ast_1$ is the norm-closed convex hull of its extreme points. Hence by the arguments similar to the ones given during the proof of Theorem~\ref{T2}, $K$ is a simplex. Now, one proceeds as in the proof of Theorem~\ref{T7}; we note that the arguments  in the latter half of the proof remain valid under the weakly Hahn-Banach smooth assumption. Thus, we get the required conclusion.
\end{proof}
The next two results further illustrate the importance of an embedding's specific nature. We recall from Section 2 that under the canonical embedding, $K$ is always a weak$^\ast$-closed face of $A(K)^\ast_1$. Additionally, a Choquet simplex whose set of extreme points is closed is called a Bauer simplex. 
\begin{thm}\label{T8}
Let $K$ be a compact convex set and suppose that $\partial_e K$ is a closed set. Consider $A(K)$ as a subspace of $C(\partial_e K)$ under the restriction embedding. Then \tfae.
\bla
\item All positive linear functionals in $S(A(K)^\ast)$ have unique positive extensions in $C(\partial_e K)^\ast$.
\item $A(K)\cong C(\partial_e K)$, i.e., the restriction map is an onto isometry.
\item $K$ is a Bauer simplex.
\el
\end{thm}
\begin{proof}
$(a)\Ra (c).$ For the proof of this part, we will use Choquet's theorem (see \cite{AE} Theorem~I.6.6). Let $x\in K$ and $\mu_1, \mu_2\in P(K)$ be maximal probability measures such that the resultants, $r(\mu_1)=r(\mu_2)=x$.  As $\partial_e K$ is a closed set and the measures are maximal, by Proposition I.4.6 in \cite{A},  we have the closed support, $Supp(\mu_i)\ci \partial_e K$, i=1, 2. Let $a\in A(K)$. Then $a(x)=\int_{\partial_e K}ad\mu_1=\int_{\partial_e K}ad\mu_2$. Thus $\mu_1, \mu_2$ are two positive extensions of $\de_x$ over $C(\partial_e K)$. From $(a)$ we get $\mu_1=\mu_2$. Hence $K$ is a simplex. Since $\partial_eK$ is a closed set, we get that $K$ is a Bauer simplex.
\vskip .5em
$(c)\Ra (b).$ Since $\partial_e K$ is compact, $A(K)$ is lattice isometric to $C(\partial_e K)$ (see \cite{AE} Theorem~II.7.5).
\vskip .5em
$(b)\Ra (a).$ This is clear.
\end{proof}
Let $\Omega$ be a compact set. Let $P(\Omega)$ be the set of probability measures in $C(\Omega)^\ast_1$. We now consider an isometric embedding $\Phi:A(K)\ra C(\Omega)$, such that $\Phi ({\bf 1})={\bf 1}$. Note that in this case $\Phi^*$ maps  $P(\Omega)$ into $K$ (as identified as states in $A(K)^\ast_1$). This is because, for any $\mu \in P(\Omega)$, $\|\Phi^\ast(\mu)\|= 1 = \Phi^\ast(\mu)({\bf{1}})=\mu(\Phi({\bf{1})})=\mu(\bf{1})= 1$. Since $\Phi$ is one-to-one, the range of  $\Phi^\ast$ is weak$^\ast$-dense in $A(K)^\ast$.  By the closed range theorem (see \cite{H}) we have, $\Phi^\ast(C(\Omega)^\ast)$ is a weak$^\ast$-closed set. Thus $\Phi^\ast$ is an onto map. In the proof of the following theorem, we again use the identification of the state space of $A(K)$ with $K$.
\vskip 1em
We further note that if $a \in A(K)_1$ and $a \geq 0$, then $0\leq a\leq {\bf1}$. Thus $0\leq {\bf 1}-a\leq {\bf 1} $ so that $\|{\bf 1}-a\|\leq 1$. Since $\Phi$ is an isometry and $\Phi({\bf1})= {\bf1}$, we have $\|\Phi(a)\|
=\|a\|$ as well as $\|{\bf 1}-\Phi(a)\|=\|{\bf1}- 
a\|$. Hence $\Phi(a) \geq 0$.
\begin{thm}\label{T9}
Let $K$ be a compact convex set and let $\Omega$ be a compact Hausdorff space. Let $\Phi: A(K)\ra C(\Omega)$ be an isometric embedding. Suppose $\Phi ({\bf 1})={\bf 1}$. If every positive linear functional on $\Phi (A(K))$ has a unique positive Hahn-Banach extension to $C(\Omega)$, then $K$ is a Bauer simplex.
\end{thm}
\begin{proof}
We first prove that $\partial_e K$ is closed in $K$. We use techniques similar to the ones used in the proof of Theorem 3.4. Let $k_\al\ra k_0$ in $K$ where $k_\al\in \partial_e K$. For $k\in K$, we continue to denote by $\delta_k$ the functional on $\Phi(A(K))$ defined by $\delta_k(\Phi(a))= a(k)$, for $a \in A(K)$. Let $\tilde{\de_{k}}$ be the unique positive extension of $\de_{k}$ to $C(\Omega)$ and from our assumption $\tilde{\de_{k}}$  is in the Bauer simplex $P(\Omega)$. As in the proof of Theorem 3.4, $\tilde{\de_{k_\al}}\ra \tilde{\de_{k_0}}$ in weak$^\ast$-topology in $C(\Omega)^*$. From Lemma~\ref{L1}, we get that $\tilde{\delta_{k_\alpha}}\in \partial_e (C(\Omega)_1^*)$, for all $\alpha$.
As $\partial_e (C(\Omega)_1^*)$ is weak$^*$-closed, we have that $\tilde{\de_{k_0}}\in\partial_e C(\Omega)_1^*$. Now $\Phi^*(\tilde{\de_{k_0}})=\de_{k_0}$. If possible assume that there exist $k_1, k_2\in K$ such that $k_0=\frac{k_1+k_2}{2}$. Then there exists $\tilde{\de_{k_i}}\in C(\Omega)_1^*$ such that $\Phi^* (\tilde{\de_{k_i}})=\de_{k_i}$, for $i=1, 2$. Thus we have $\frac{\tilde{\de_{k_1}}+\tilde{\de_{k_2}}}{2}=\tilde{\de_{k_0}}$. This contradicts the fact that $\tilde{\de_{k_0}}$ is an extreme point of $C(\Omega)_1^*$. Hence $k_0\in\partial_e K$.
It remains to prove that $K$ is a simplex.
\vskip .5em
{\sc Claim:~} Every closed face of $K$ is a parallel face. The proof proceeds in several steps.
\vskip .5em
Let $F$ be a closed face of $K$. For each $z\in F$, let $\widetilde{\de_z}\in C(\Omega)^*$ be as stated above. Let $\widetilde{F}=\{\widetilde{\de_z}:z\in F\}$, we will show that $\widetilde{F}$ is a face of $P(\Omega)$. Note that $\Phi^*\widetilde{\de_z}(a)=\widetilde{\de_z}(\Phi (a))$ for $a\in A(K)$. Hence we get $\Phi^*(\widetilde{\de_z})= \delta_z$. Since $A(K)$ separates points of $K$, we get $\Phi^* (\widetilde{F})=F$. Suppose $\Theta_i\in P(\Omega)$, $i=1, 2$ such that $\frac{\Theta_1+\Theta_2}{2}\in \widetilde{F}$. Then it is clear that $\Phi^* (\Theta_i)\in A(K)_1^*$ are functionals which satisfy $\frac{\Phi^*(\Theta_1)+\Phi^*(\Theta_2)}{2}\in \Phi^* (\widetilde{F})=F$. Since $F$ is also a face of $A(K)^\ast_1$, it follows that $\Phi^* (\Theta_1), \Phi^* (\Theta_2)\in F$. As $\Theta_i$ exists uniquely for which $\Phi^*(\Theta_i)= \delta_{s_i}$, for some $s_i \in F$, for $i=1,2$, we have $\Theta_i=\widetilde{\de_{s_i}}$, for $i=1, 2$. Consequently, $\Theta_i\in \widetilde{F}, i=1,2$ and this proves that $\widetilde{F}$ is a face, and arguments similar to the ones given before also yield that it is a closed face. Now identifying $C(\Omega)=A(P(\Omega))$, from our concluding remarks in Section 2, $\widetilde{F}$ is a parallel face of $P(\Omega)$.
\vskip .5em

Let $G\ci P(\Omega)$ be a face for which $P(\Omega)=co \{G\cup \widetilde{F}\}$. Consider $\Phi^\ast(G) \ci A(K)^\ast_1$.
\vskip .5 em
 The following steps yield our {\sc Claim}.
{\sc Step 1:~} $\Phi^* (G)$ is a face of $K$.
{\sc Step 2:~} $co \{F\cup \Phi^*(G)\}= K$.
{\sc Step 3:~} $F$ is a parallel face of $K$.
\vskip .5em
Proof of Step 1: Let $\La_i\in A(K)_1^*$ for $i=1, 2$ be such that $\frac{\La_1+\La_2}{2}\in \Phi^* (G)$. Now $\La_i\in K\cup (-K)$, $i=1, 2$. Since $\La_i$'s have unique extensions to $C(\Omega)$, there exists $\widetilde{\La_i}\in P(\Omega)\cup -P(\Omega)$ such that $\Phi^* (\widetilde{\La_i})=\La_i$, $i=1,2$. Hence we have $\frac{\Phi^*(\widetilde{\La_1})+\Phi^*(\widetilde{\La_2})}{2}\in \Phi^* (G)$. It is clear that $\Phi^*$ is one-to-one on the unique extensions. This concludes that $\frac{\widetilde{\La_1}+\widetilde{\La_2}}{2}\in G$. Since $G$ is a face, $\widetilde{\La_1},\widetilde{\La_2} \in G$. Hence, the {\sc Step 1} follows.

Proof of Step 2: Let $k\in K$. By our assumption, there exists unique Hahn-Banach extension $\La\in P(\Omega)$ such that $\Phi^* (\La)=\de_k$. Since $\widetilde{F}$ is a parallel face $\de_k \in co\{F \cup \Phi^\ast(G)\}$. Hence  {\sc Step 2} follows.
\vskip .5em

Proof of Step 3:
Now suppose for $k\in K$ there exist distinct representations $k=\la v_1+(1-\la) v_2=\mu u_1+(1-\mu) u_2$, for $0\leq \la, \mu\leq 1$ and $\la\neq\mu$, where $v_1, u_1\in F$ and $v_2, u_2\in \Phi^* (G)$. Then consider the corresponding extensions of $\de_{v_i}, \de_{u_i}$, $i=1, 2$ to $C(\Omega)$. As both the functionals $\la\widetilde{\de_{v_1}}+(1-\la)\widetilde{\de_{v_2}}$ and $\mu\widetilde{\de_{u_1}}+(1-\mu)\widetilde{\de_{u_2}}$ are in $C(\Omega)_1^*$ represent the extensions of $\de_k$ to $C(\Omega)$, denoted as $\widetilde{\de_k}$. Since $\la\neq\mu$, this contradicts the fact that $\widetilde{F}$ is a parallel face of $P(\Omega)$.
Hence, the {\sc Step 3} follows.
This completes the proof of the {\sc Claim}. As every closed face is a parallel face, it follows from  Theorem III. 2. 5 in \cite{AE}, that $K$ is a simplex.
\end{proof}
 Operator-theoretic versions, where one studies unique, positive, and norm-preserving extensions of positive operators, of the results considered here are still open. The following is a folklore result. Let $K$ be an infinite-dimensional Choquet simplex, whose extreme boundary is not closed. Consider the identity map $i: A(K) \rightarrow A(K)$. There is no linear operator $T: C(K) \rightarrow A(K)$ such that $\|T\|=1$ and $T= i$ on $A(K)$. To see this, note that the existence of a norm preserving or equivalently a linear positive extension implies that $A(K)$ is the range of a contractive projection on $C(K)$. It follows from Theorem 5 in Section 10 of \cite{L} (see also the classification scheme in \cite{Li}), since $A(K)$ is the range of a contractive projection in a space of continuous functions, $\partial_e A(K)^\ast_1 \cup \{0\}$ is a weak$^\ast$-compact set. Since ${\bf 1} \in A(K)$, $0$ cannot be an accumulation point. Now it is easy to see that $\partial_e K$ is closed, yielding the required contradiction.
\begin{prob}
 In view of the operator versions of unique isometric extensions between spaces of compact operators and spaces of operators,  from \cite{Rao1}, it would be interesting to find conditions under which an order-preserving isometric embedding $\Phi: A(K) \rightarrow C(K)$ has a unique extension to an order-preserving isometric embedding from $C(K)$ to $C(K)$?
	
\end{prob}
\vskip 2em
The second author's work is supported by a 3-year project, `Classification of Banach spaces using differentiability,' funded by the Anusandhan National Research Foundation under Core Research Grant CRG-2023-000595. He also thanks the Department of Mathematics, IITH, for their hospitality when part of this work was done. The authors thank the referees for their comments and suggestions.

\vskip 1em
\begin{center}
  \textbf{DECLARATION}
\end{center}
The authors declare that no data were generated during this investigation, and hence, the 'data sharing' is not applicable for this research. In addition, the authors declare no conflicts of interest.

\end{document}